\documentclass[12pt]{amsart}
\usepackage{amsfonts,latexsym,amsthm,amssymb}
\usepackage[all]{xy}

\input epsf
\CompileMatrices

\newtheorem{theorem}{Theorem}[section]

\newtheorem{prop}[theorem]{Proposition}

{\theoremstyle{definition} }
{\theoremstyle{remark} 
\newtheorem{example}[theorem]{Example}}

\newcommand{\Abb}{{\mathbb{A}}}
\newcommand{\Cbb}{{\mathbb{C}}}
\newcommand{\Qbb}{{\mathbb{Q}}}

\newcommand{\Pbb}{{\mathbb{P}}}
\newcommand{\Zbb}{{\mathbb{Z}}}

\newcommand{\one}{1\hskip-3.5pt1}
\newcommand{\csm}{c_{\text{\rm SM}}}
\newcommand{\cf}{c_{\text{\rm F}}}

\title{Computing characteristic classes of projective schemes}
\author{Paolo Aluffi}
\address{Max-Planck-Institut f\"ur Mathematik, Bonn, Germany}
\address{Florida State University, Tallahassee, Florida, USA}

\makeindex

\begin{document}

\begin{abstract}
We discuss an algorithm computing the push-forward to projective space
of several classes associated to a (possibly singular, reducible,
nonreduced) projective scheme. For example, the algorithm yields the
topological Euler characteristic of the support of a projective scheme
$S$, given the homogeneous ideal of $S$. The algorithm has been
implemented in {\tt Macaulay2.\/}
\end{abstract}

\maketitle

\section{Introduction}

\subsection{}
In this article we describe an algorithm computing, among other
things, the topological Euler characteristic of the support of a
projective scheme $S$ over $\Cbb$. In fact, we will compute the
push-forward to $\Pbb^n$ of the {\em Chern-Schwartz-MacPherson\/}
class $\csm(S)$ of the support of $S$, given the ideal $I$ of $S$ in
$\Pbb^n$; as is well known, the Euler characteristic equals the degree
of the component of dimension~0 of $\csm(S)$. We also include a
computation of the push-forward of the {\em (Chern-)Fulton\/} class
$\cf(S)$ of $S$; when $S$ is nonsingular, this provides a different
way to compute the Euler characteristic of~$S$.

Other algorithms computing the Euler characteristic of a (possibly
singular) scheme are somewhat indirect (see Uli Walther's contribution
to \cite{compalggeom}, as well as \cite{MR1874280}). The nonsingular
case can be treated by computing the Hodge numbers $h^{ij}$. Even in
the nonsingular case, however, we are not aware of algorithms yielding
(the degrees of) the Chern classes of $S$; for a nonsingular variety,
the outputs of our algorithms for $\csm(S)$ and $\cf(S)$ coincide, and
consist precisely of this information.

\subsection{}
The main ingredients to our algorithms are the results of
\cite{MR2001i:14009} and \cite{incexcI}, and explicit computations of
{\em Segre classes.\/} The considerations in \cite{incexcI} reduce the
problem of the computation of $\csm(S)$, for $S\subset \Pbb^n$, to the
case in which $S$ is a hypersurface in $\Pbb^n$; the main result of
\cite{MR2001i:14009} translates this case to the computation of a
Segre class; and a close look at Segre classes in $\Pbb^n$ reveals
that tools such as Macaulay2 (\cite{M2}) are capable of computing them.

In fact the ability to compute Segre classes appears to us of independent 
interest, for example in view of potential applications to enumerative 
geometry. An immediate application to characteristic classes yields 
the Fulton class $\cf(S)$ of $S$ (term by which we refer to the class 
introduced by William Fulton in \cite{MR85k:14004}, Example~4.2.6(a)).

\subsection{}
The article is organized as follows. In \S\ref{csm} we describe the
algorithm computing $\csm(S)$, and hence $\chi(S)$. We have given this
discussion a prominent place since it may be the item of more
immediate interest in the paper; but in fact at one key step in the
proof of the main result in \S\ref{csm}, and in the resulting
algorithm, we will borrow some material from the following
\S\ref{segre}. The algorithm is summarized in \S\ref{summary}. We end
\S\ref{csm} by pointing out that judicious use of the algorithm yields
the computation of Euler characteristics of {\em affine\/} schemes
(over a field) as well---and hence in principle of arbitrary schemes,
as every scheme is the disjoint union of affine ones.

In \S\ref{segre} we discuss the problem of computing more general Segre 
classes. Serious applications are so far severely limited by technological 
constraints. However, one subproduct of the discussion in \S\ref{segre} is 
the algorithm giving Fulton class. 

Several concrete examples are given in \S\ref{examples}. Among these, we 
mention the computation of {\em Milnor classes\/} of a projective scheme, as 
these have been the subject of rather intense work in recent years.
Briefly, the Milnor class measures the difference between 
Chern-Schwartz-MacPherson and Fulton classes of a singular variety. For 
complete intersections, Shoji Yokura (\cite{MR2000e:14011}) 
has identified the computation of these classes as a Verdier-Riemann-Roch 
type problem. The most general results obtained in this direction are in 
the recent \cite{math.AG/0202175}; for surveys of work on Milnor
classes, see \cite{MR2001k:14014} and \cite{MR2002c:14001}.

\subsection{}
We have implemented the algorithms described in this paper in
Macaulay2. Our code (and, we hope, future improvements) is available at

{\tt http://www.math.fsu.edu/\char126aluffi/CSM/CSM.html}

\noindent In any case, the reader should have no difficulties
translating the discussion presented in this paper into working
routines in Macaulay2 or other commutative algebra/algebraic geometry
symbolic packages.\vskip 6pt

\noindent{\bf Acknowledgments} I thank the Max-Planck-Institut f\"ur
Mathematik in Bonn, Germany, for the hospitality and support, and
Florida State University for granting a sabbatical leave in 2001-2.

\section{Chern-Schwartz-MacPherson classes and the Euler 
characteristic}\label{csm}

\subsection{}\label{ground}
Throughout the paper, $i:S \hookrightarrow \Pbb^n=\Pbb^nk$ will denote a 
closed embedding; in this section $k$ will be a field of characteristic~0. 
We will let $I=k[z_0,\dots,z_n]$ be a homogeneous ideal defining $S$.

The output of our computations will be classes in the Chow group of 
$\Pbb^n$. Denoting by $H$ the hyperplane class, this is 
$\Zbb[H]/(H^{n+1})$: classes in $\Pbb^n$ will be written as polynomials 
of degree $\le n$ in $H$, with integer coefficients:
$$a_0+a_1H+\dots+a_n H^n\quad.$$
The {\em degree\/} of a class, denoted $\int$, will be the coefficient of 
$H^n$ in such an expression.

\subsection{}
If $S$ is a nonsingular variety, we may consider its (total, homology) 
{\em Chern class\/} $c(TS)\cap [S]$. {\em Chern-Schwartz-MacPherson\/}
classes provide a notion agreeing with $c(TS)\cap [S]$ when $S$ is
nonsingular, but defined regardless of the nonsingularity of
$S$. Further, they satisfy a clever functorial prescription, which we
quickly summarize. 

Denote by $\csm(S)$ the Chern-Schwartz-MacPherson class of $S$, and extend 
this definition to {\em constructible functions\/} by setting
$$\csm(\sum_{V\subset S} m_V \one_V)=\sum_V m_V \csm(V)\quad.$$
Here the sum is finite, $V$ are closed subvarieties of $S$, $m_V\in 
\Zbb$, and $\one_V$ denotes the function that is $1$ along $V$ and $0$ 
outside of $V$. This defines a homomorphism of abelian groups $\mathcal 
C(S) \to \mathcal A(S)$ for every $S$, where $\mathcal C$, $\mathcal A$ 
denote respectively the functor of constructible functions (with 
push-forward defined by Euler characteristic of fibers) and the Chow group 
functor. But in fact
$$\csm : \mathcal C \leadsto \mathcal A$$
defines a natural transformation: this was proved by Robert MacPherson in 
the article where the classes are introduced. For MacPherson's construction 
of $\csm$, and for more information, we address the reader to the original
\cite{MR50:13587}, or to \cite{MR91h:14010} (extending the theory to 
arbitrary algebraically closed field of characteristic~0); and to
\cite{MR2002c:14001} for a comparison with the different approach of 
Marie-H\'el\`ene Schwartz, in fact predating MacPherson's work. Regardless of 
the approach, at the moment the theory of Chern-Schwartz-MacPherson classes 
has only been studied in characteristic~$0$, and this is why we assume that 
our ground field is of characteristic~$0$ in this section.

\subsection{}
In fact, the theory is usually only applied to {\em reduced\/} schemes.
More generally, we take $\csm(S)$ to be the Chern-Schwartz-MacPherson 
class of the support $S_{\text{red}}$ of $S$. As a very particular case of the 
functoriality of Chern-Schwartz-MacPherson classes, consider the constant map 
on a proper scheme $S$,
$$\kappa:S \to \text{point}\quad.$$
Then the covariance of $\csm$ for $\kappa$ amounts to
$$\kappa_*\csm(S)=\csm(\kappa_*\one_S)=\csm(\chi(S_{\text{red}})
\one_{\text{point}})=\chi(S_{\text{red}}) [\text{point}]\quad,$$
and in particular
$$\int \csm(S)=\chi(S_{\text{red}})\quad.$$
With $S$ projective, and $i:S\to \Pbb^n$ a closed embedding, this says
that
$$\chi(S_{\text{red}})=\int i_*\csm(S)\quad:$$
that is, the topological Euler characteristic of the support of $S$
equals the coefficient of $H^n$ in
$$i_*\csm(S)=c_0+c_1 H+\dots+c_n H^n\quad.$$
Computing this class is our main goal.

We note in passing that the computations we will describe can all be
performed over any field over which $S$ is defined. Thus, a tool such
as Macaulay2 will be able to compute the topological Euler
characteristic of a scheme $S\subset\Pbb^n\Cbb$ by working over $\Qbb$
(for example), so long as $S$ is in fact defined over $\Qbb$.

\subsection{}
We will now describe a procedure computing $i_*\csm(S)$, given a 
homogeneous ideal 
$$I=(F_1,\dots,F_r)$$ 
defining $S$ in $\Pbb^n$. Write
$$S=X_1\cap\cdots\cap X_r$$
where $X_i$ is the hypersurface defined by $F_i$. Applying
`inclusion-exclusion', it suffices then (see \cite{incexcI}) to
compute $i_*\csm(X)$ for $X$ ranging over the unions
$X_{i_1}\cup\cdots\cup X_{i_s}$, $1\le s\le r$. In other words, the
problem of computing $i_*\csm(S)$ is readily reduced to the computation
of $i_*\csm(X)$ for $X$ a {\em hypersurface\/} in $\Pbb^n$.

If $X$ is a hypersurface, the computation of $\csm(X)$ is reduced to the 
computation of the Segre class of the {\em singularity subscheme\/} of $X$ by 
the main result of \cite{MR2001i:14009}. The form taken by this result when 
the ambient nonsingular variety is projective space is the following.

\subsection{}
Let $X\subset \Pbb^n$ be a hypersurface, with homogeneous ideal
$(F)\subset k[z_0,\dots,z_n]$. Consider the rational map
$$\Pbb^n \dashrightarrow \Pbb^N=\Pbb^n$$
defined by 
$$p \mapsto \left( \tfrac{\partial F}{\partial z_0}_{|_p}: \dots : 
\tfrac{\partial F}{\partial z_n}_{|_p}\right)\quad.$$
We let $\Gamma$ be the (closure of the) graph of this map. Viewing 
$\Pbb^n\times \Pbb^N$ as a $\Pbb^N$-bundle over $\Pbb^n$, we are interested 
in what we elsewhere call the {\em shadow\/} of the class of $\Gamma$: 
that is, letting $K$ be the pull-back of the hyperplane class from 
$\Pbb^N$, the class
$$G=g_0+g_1 H+\dots +g_n H^n$$
in $A_*\Pbb^n$, where $g_i$ is the degree of the image in $\Pbb^n$ of 
$K^i\cdot [\Gamma]$.

\begin{theorem}\label{csmthm}
With the notations introduced above,
$$i_*\csm(X)=(1+H)^{n+1}-\sum_{i=0}^n g_d\, (-H)^d (1+H)^{n-d}\quad.$$
\end{theorem}

\begin{proof}
By Theorem~I.4 in \cite{MR2001i:14009},
$$i_*\csm(X)=c(T\Pbb^n)\cap i_*\left(s(X,\Pbb^n)+c(\mathcal O(X))^{-1}\cap
\left( s(Y,\Pbb^n)^\vee\otimes \mathcal O(X)\right)\right)$$
where $Y$ denotes the singularity subscheme of $X$; this is the scheme 
defined by the vanishing of the partials of $F$. By Proposition~\ref{redux} 
in the next section, $i_*s(Y,\Pbb^n)$ can be recovered from the class
$G=g_0+g_1 H+\dots +g_n H^n$:
$$i_* s(Y,\Pbb^n)=1-c(\mathcal O(dH))^{-1}\cap\left( G\otimes \mathcal
O(dH) \right)\quad,$$
where $d=\deg X-1$ (so $\mathcal O(X)=\mathcal O((d+1)H)$).

The manipulations massaging this formula into the one given in the statement 
are streamlined by using Proposition~1 in \cite{MR96d:14004}:
$$i_* s(Y,\Pbb^n)^\vee=1-c(\mathcal O(-dH))^{-1}\cap\left( G^\vee\otimes 
\mathcal O(-dH) \right)$$
$$i_* s(Y,\Pbb^n)^\vee\otimes\mathcal O(X)=1-\frac{c(\mathcal O(X))} 
{c(\mathcal O(H))}\cap \left(G^\vee\otimes \mathcal O(H)\right)$$
$$c(\mathcal O(X))^{-1}\cap \left(i_* s(Y,\Pbb^n)^\vee\otimes\mathcal O(X) 
\right)=c(\mathcal O(X))^{-1} - c(\mathcal O(H))^{-1}\left(G^\vee\otimes
\mathcal O(H)\right)\quad,$$
and hence
$$\aligned
i_*\csm(X)& =(1+H)^{n+1}\left(c(\mathcal O(X))^{-1}\cap[X]+
c(\mathcal O(X))^{-1} - c(\mathcal O(H))^{-1}\left(G^\vee\otimes
\mathcal O(H)\right)\right)\\
&=(1+H)^{n+1}-(1+H)^n (G^\vee\otimes\mathcal O(H))
\endaligned$$
which translates into the formula given in the statement.
\end{proof}

\subsection{}\label{summary}
Therefore, up to the bookkeeping of inclusion-exclusion and to trivial 
algebraic manipulations, the problem of computing $i_*\csm(S)$ is reduced 
by Theorem~\ref{csmthm} to the computation of the shadow $G$ of the graph 
$\Gamma$ of a rational map. This is the key ingredient, and since it 
yields more generally the Segre class of any closed subscheme of
$\Pbb^n$ we discuss it separately, in \S\ref{segre}. 

Summarizing: given the ideal $I=(F_1,\dots,F_r)$ of $S$, an algorithm 
computing $i_*\csm(S)$ will
\begin{itemize}
\item list all products $F=\prod_{i_1<\dots<i_s} F_{i_1}\cdot \cdots \cdot 
F_{i_s}$;
\item for each such $F$, compute the jacobian ideal $J=(\frac{\partial F}
{\partial z_0},\dots, \frac{\partial F}{\partial z_n})$, and apply the 
procedure described in \S\ref{segre} in order to compute the corresponding 
class $G$;
\item apply Theorem~\ref{csmthm} to this class, and obtain $i_*\csm(X)$ 
for the hypersurface $X$ corresponding to $F$;
\item apply inclusion-exclusion to reconstruct $i_*\csm(S)$.
\end{itemize}
The coefficient of $H^n$ in $i_*\csm(S)$ gives the Euler characteristic of 
the support of $S$.

\subsection{}
If $r$ is the number of generators of the ideal of $S$, one `shadow'
computation is required for each of the $2^r-1$ hypersurfaces invoked by 
inclusion-exclusion. This causes an exponential slow-down of the procedure 
as the codimension of $S$ increases. 

It is somewhat amusing that the result of the computation, that is,
$i_*\csm(S)$, only depends on the support of $S$, even if in no place
does the algorithm explicitly compute the support of $S$, or of the
hypersurfaces $X$ considered at intermediate stages. In fact,
introducing intermediate computations of supports may speed up the
algorithm: any procedure `simplifying' the input $I$---in the sense of
reducing the number and degree of the generators, without altering
the radical of $I$---should lead to an increase in the efficiency of the
procedure.

\subsection{}\label{affine}
The procedure is easily adapted to the computation of the Euler 
characteristic of (the support of) a closed subscheme $S$ of 
affine space $\Abb^n$, given its defining ideal. This can be done in 
several ways: for example, one may homogenize the ideal of $S$, obtaining 
the closure $\overline S\subset\Pbb^n$; then multiply this ideal by the 
equation of the hyperplane $L$ at infinity, obtaining the union 
$\overline S\cup L\subset \Pbb^n$; and then compute
$$\chi(S)=\chi(\overline S\cup L)-(n+1)\quad.$$

As an alternative, one may intersect with the hyperplane at infinity, 
obtaining a `limit' subscheme $\underline S\subset\Pbb^{n-1}$; and then
$$\chi(S)=\chi(\overline S)-\chi(\underline S)\quad.$$
This approach appears to be much faster in practice.

\section{Computing Segre classes of subschemes of $\Pbb^n$}\label{segre}

\subsection{}
We can now lift the restriction on the characteristic of the ground field 
$k$, as they are irrelevant for the considerations in this
section. Again we let $i:S\hookrightarrow \Pbb^n$ be a closed
embedding of a scheme $S$ in projective space $\Pbb^n=\Pbb^nk$; our
goal is to give an explicit procedure computing the push-forward 
$$i_* s(S,\Pbb^n)\in A_*\Pbb^n$$
of the Segre class $s(S,\Pbb^n)$ of $S$ in $\Pbb^n$. By 
Proposition~\ref{redux} this will be reduced to the computation of a 
`shadow', as has been the case in \S\ref{csm}; we will then discuss the 
computation of shadows, in \S\ref{step1} and ff.

\subsection{}
Let $I=(f_0,\dots,f_N)\subset k[z_0,\dots,z_n]$ be a homogeneous ideal
defining $S$. We may and will assume that the generators $f_i$ are
all of the same degree $r$; in other words, we write $S$ as the
zero-scheme of a section of $\mathcal O(rH)^{\oplus(n+1)}$:
$$(f_0,\dots,f_N):\mathcal O_{\Pbb^n} \rightarrow \mathcal
O(rH)^{\oplus(n+1)}\quad.$$
Projectivizing, we get a rational map
$$\Pbb^n \dashrightarrow \Pbb^N$$
and we let 
$$\Gamma_I\subset \Pbb^n\times \Pbb^N$$ 
denote the (closure of the) graph of this map. Denote by $K$ the
pull-back of the hyperplane class from the $\Pbb^N$ factor, and by
$\pi$ the projection $\Gamma_I \to \Pbb^n$. The shadow of $\Gamma_I$ is 
the class
$$G=g_0+g_1 H+\dots+g_n H^n\in A_*\Pbb^n\quad,$$
where $g_i$ is the degree of $\pi_*(K^i\cdot [\Gamma_i])$. 

\subsection{}
Now we can state and prove the simple result used in the proof of 
Theorem~\ref{csmthm}. The statement, as that proof, uses the 
notations from \cite{MR2001i:14009}.

\begin{prop}\label{redux}
With the above positions,
$$i_* s(S,\Pbb^n)=1-c(\mathcal O(rH))^{-1}\cap\left( G\otimes \mathcal
O(rH) \right)\quad.$$
\end{prop}

\begin{proof}
By construction, the graph $\Gamma_I$ is isomorphic to the blow-up of
$\Pbb^n$ along $S$, and the class of the exceptional divisor $E$ on
$\Gamma_I$ equals the restriction of $c_1(\mathcal O(-1))$ from
$\Pbb(\mathcal O(rH)^{\oplus(n+1)})\cong \Pbb^n$. Chasing this
identification, we see that the class of $E$ is $rH-K$. Hence using
\cite{MR85k:14004}, Corollary~4.2.2:
$$s(S,\Pbb^n)=\pi_*\frac {[E]}{1+E}=\pi_*\frac{[rH-K]}{1+rH-K}\quad.$$
Pushing forward to $\Pbb^n$, this can be manipulated as follows:
$$\aligned
\pi_*\frac{[rH-K]}{1+rH-K} &=\pi_*\left([\Gamma_I]-\frac
1{1+rH-K}\cdot [\Gamma_I]\right)\\
&=1-\pi_*\left(\frac 1{1+rH}\cdot \frac {1+rH}{1+rH-K}\cdot
[\Gamma_I] \right)\\
&=1-c(\mathcal O(rH))^{-1}\pi_* \cap\left(\left( \frac 1{1-K}\cdot
[\Gamma_I]\right)\otimes \mathcal O(rH)\right)\\
&=1-c(\mathcal O(rH))^{-1}\cdot \left(G\otimes \mathcal O(rH)
\right)
\endaligned$$
as claimed.
\end{proof}

\subsection{}
The upshot of Theorem~\ref{csmthm} and Proposition~\ref{redux} is that 
we can compute Chern-Schwartz-MacPherson classes and Segre classes (and 
hence Fulton classes) if we can extract the integers $g_i$ giving the 
coefficients of the class $G$ determined by a graph $\Gamma_I$. That is, we 
must be able to 
\begin{itemize} 
\item obtain $\Gamma_I$ explicitly;
\item intersect $\Gamma_I$ with general hyperplanes;
\item project the intersections down to $\Pbb^n$;
\end{itemize}
and compute the degree of these projections.

Each of these steps is easily implemented in any of the standard
symbolic computations packages; we briefly discuss this in the following 
subsections.

\subsection{Obtaining $\Gamma_I$ explicitly}\label{step1}
A bihomogeneous ideal for the graph $\Gamma_I$ can be given in the
ring 
$$k[t_0,\dots,t_N,z_0,\dots,z_n]$$
by the following trick (going back at least as far as
\cite{MR31:1275}, proof of Lemma~1): adjoin an auxiliary variable $u$
to the ring, and consider the ideal
$$J=(t_0-uf_0,\dots,t_N-uf_N)$$
in the extended ring. Then the ideal for $\Gamma_I$ is the contraction
$$J_0:=J\cap k[t_0,\dots,t_N,z_0,\dots,z_n]\quad.$$
Indeed, $J$ is the kernel of the homomorphism
$$k[u,t_0,\dots,t_N,z_0,\dots,z_n]\rightarrow k[u,z_0,\dots,z_n]$$
obtained by mapping $t_i$ to $u f_i$; a polynomial $P\in 
k[t_0,\dots,t_N,z_0,\dots,z_n]$ maps to $0$ by this map if and only if
$P$ vanishes whenever $(t_0:\dots:t_N)=(f_0:\dots:f_N)$.

The ideal $J_0$ of $\Gamma_I$ can thus be obtained by standard
elimination theory: choose a monomial order so that $u$ precedes the
other variables; compute a Gr\"obner basis for $J$; and eliminate $u$
to obtain the intersection of $J$ with the ring
$k[t_0,\dots,t_N,z_0,\dots,z_n]$.

Needless to say, this operation is rather computationally expensive. 
Of course, any other algorithm computing the Rees algebra of $I$ can
be employed here; the topic is treated extensively in
\cite{MR99c:13048}, \S7.2.

\subsection{Intersecting $\Gamma_I$ with general hyperplanes}\label{step2}
Programs such as Macaulay2 include the option of producing `random'
elements of given degree in a ring; for $i=1,\dots,n$ we can inductively
set
$$J_i:=\text{saturate}(J_{i-1}+(\ell_i),(t_0,\dots,t_N))
\quad,$$
where $\ell_i=\ell_i(t_0,\dots,t_N)$ is a random linear polynomial in
$k[t_0,\dots,t_N]$, and the saturation is necessary to remove possible
components in the intersection supported on the irrelevant ideal,
see below.

Of course we have to take care that random is sufficiently random. For
the purposes of this computation, a hyperplane is general if it
does not contain any component of the object it is intersecting, that
is, if it is not contained in any of the associated primes of the
corresponding ideal. This can be explicitly checked, for example by
making sure that the dimension decreases upon intersecting with the
hyperplane. Thus the ideals $J_i$ can be obtained as above, by
producing enough random $\ell_i$ until a general one is found.

As for the saturation, the manipulation of the ideals in
$k[t_0,\dots,t_N,z_0,\dots,z_n]$ amounts to working in
$\Abb^{n+1}\times \Abb^{N+1}$. Saturating with respect to the
irrelevant ideal $(t_0,\dots,t_N)$ guarantees that there is a
bijection between the components of the subscheme defined by $J_i$ in
$\Abb^{n+1}\times \Abb^{N+1}$ and those (about which we are
interested) in $\Abb^{n+1}\times \Pbb^N$.

\begin{example}\label{satur}
Here is an example showing that extra components may indeed
appear. Consider $I=(z_0,z_1)$ in $k[z_0,z_1]$. Then, with
notations as above,
$$J_0=(z_0t_1-z_1t_0)\quad.$$
Intersecting by $t_0$ does decrease the dimension (so that $t_0$
is general in the above sense), but creates a component supported on
the irrelevant ideal:
$$(z_0t_1-z_1t_0)+(t_0)=(t_0,t_1)\cap (z_0,t_0)\quad.$$
Saturating with respect to $(t_0,t_1)$ eliminates such spurious
components.
\end{example}

\subsection{Projecting down to $\Pbb^n$}\label{step3} This is also 
done by elimination theory. Once $J_i$ is obtained, we can ask for 
the Gr\"obner basis with respect to a monomial ordering in which
$t_0,\dots,t_N$ precede $z_0,\dots,z_n$, then eliminate
$t_0,\dots,t_N$. This computes
$$J_i\cap k[z_0,\dots,z_n]\quad,$$
that is, the homogeneous ideal in $\Pbb^n$ of the projection of the
$i$-th linear section.

\subsection{} 
Programs such as Macaulay2 compute the degree of the scheme defined by a 
given homogeneous ideal without difficulty. Applying this to the ideal 
obtained in the previous step produces the list of integers
$$g_0=1\,,\, g_1\,,\, \dots\,,\, g_n$$
needed in Theorem~\ref{csmthm} and Proposition~\ref{redux}.

\subsection{}\label{fulton}
Every scheme $S$ embeddable in a nonsingular variety $M$ has an intrinsic {\em 
Fulton class\/} 
$$\cf(S)=c(TM)\cap s(S,M)$$
(see \cite{MR85k:14004}, Example 4.2.6(a)). As $i_* s(S,\Pbb^n)$ is 
available via the procedure described above, so is
$$i_* \cf(S)=(1+H)^{n+1}\cdot i_* s(S,\Pbb^n)$$
for a projective scheme.

That $\cf(S)$ is {\em intrinsic\/} means that it does not depend on the 
chosen embedding. For example, if $S$ itself is nonsingular, then
$$\cf(S)=c(TS)\cap [S]\quad,$$
and in particular
$$\int i_*\cf(S)=\int c(TS)\cap [S]=\chi(S)$$
computes (in characteristic~0) the Euler characteristic of $S$. It
seems, however, that the computation of the Euler characteristic via
$h^{ij}$ would be much faster in this case.

Not much is known about $\cf(S)$ in general, even regarding
$\int\cf(S)$ (cf.~\cite{MR85k:14004}, Example 4.2.6(b)). If $S$ is a
local complete intersection, then $\cf(S)$ equals the class of the
virtual tangent bundle of $S$. In this case, identifying the
difference between $\cf(S)$ and the functorial $\csm(S)$ has been
identified by Yokura as a Verdier-type Riemann-Roch problem;
cf.~Example~ref{milnor}.

\section{Examples}\label{examples}

We won't reproduce here the Macaulay2 code implementing the above
steps, as further details seem unnecessary, and our code is certainly
much less than optimal. A documented copy of the code (and of future
improvements) is available at

{\tt http://www.math.fsu.edu/\char126aluffi/CSM/CSM.html}

\noindent In the present version, loading the code (named {\tt
CSM.m2}) produces several functions:
\begin{itemize}
\item {\tt segre}
\item {\tt cf}
\item {\tt csm}
\item {\tt euleraffine}
\end{itemize}
with hopefully evident meaning. The first three items accept a homogeneous 
ideal in a polynomial ring as argument; {\tt euleraffine} accepts a (not 
necessarily homogeneous) ideal in a polynomial ring.

The simple examples which follow are meant to illustrate the use of these
functions.

\begin{verbatim}
Macaulay 2, version 0.9
--Copyright 1993-2001, D. R. Grayson and M. E. Stillman
--Singular-Factory 1.3b, copyright 1993-2001, G.-M. Greuel, et al.
--Singular-Libfac 0.3.2, copyright 1996-2001, M. Messollen

i1 : load "CSM.m2"
--loaded CSM.m2
\end{verbatim}
Most of the examples are chosen in projective spaces of dimension~2,3,
and~4 over $\Bbb Q$:
\begin{verbatim}
i2 : ringP2=QQ[x,y,z]; ringP3=QQ[x,y,z,w]; ringP4=QQ[x,y,z,w,t];
\end{verbatim}
\begin{example}[Three concurrent lines in $\Pbb^3$]
The Segre class of the reduced scheme $S$ supported on three general distinct 
lines through a point in $\Pbb^3$ is computed by
\begin{verbatim}
i5 : use ringP3; segre ideal(x*y,x*z,y*z)
                   3     2
Segre class : - 10H  + 3H
\end{verbatim}
The output is written in the Chow ring of $\Pbb^3$, where $H$ denotes
the hyperplane class; thus the result is
$$i_* s(S,\Pbb^3)=3[\Pbb^1]-10[\Pbb^0]\quad.$$

The class changes if the lines become coplanar. For instance, consider the 
ideal $(z, xy(x+y))$ (in order to compute the Segre class 
in this case, the routine modifies the ideal so that all generators
have the same degree:
$$(x^2z,y^2z,z^3,z^2w,xy(x+y))\quad;$$
this ideal defines the same scheme, so it yields the same Segre class).
\begin{verbatim}
i7 : segre ideal(z,x*y*(x+y))
                   3     2
Segre class : - 12H  + 3H
\end{verbatim}
Or we may argue that since three coplanar lines form a plane curve of
degree~3, the Fulton class of $S$ must equal the class for a
nonsingular plane cubic; then use that Fulton classes are intrinsic
(see \S\ref{fulton}) to compute the Segre class in $\Pbb^3$. This
gives the same result:
$$\frac{3H^2}{(1+H)^4}=3H^2-12 H^3\quad.$$

In order to compute directly the Fulton classes for these two examples:
\begin{verbatim}
i8 : CF ideal(x*y,x*z,y*z)
                 3     2
Fulton class : 2H  + 3H

i9 : CF ideal(z,x*y*(x+y))
                 2
Fulton class : 3H
\end{verbatim}
while the Chern-Schwartz-MacPherson classes are:
\begin{verbatim}
i10 : CSM ideal(x*y,x*z,y*z)
                                    3     2
Chern-Schwartz-MacPherson class : 4H  + 3H

i11 : CSM ideal(z,x*y*(x+y))
                                    3     2
Chern-Schwartz-MacPherson class : 4H  + 3H
\end{verbatim}
This example illustrates that Chern-Schwartz-MacPherson classes are,
to some extent, `combinatorial objects': unlike Fulton classes, they
do not tell the  difference between the two configurations.
\end{example}

\begin{example}[plane cubics]
However, Fulton classes cannot tell the difference between a nonsingular 
plane cubic and a singular one. We switch to dimension 2, which speeds up 
the computations somewhat; the Fulton classes for $(x^3+y^3+z^3)$ and
$(xy(x+y))$ agree:
\begin{verbatim}
i12 : use ringP2; CF ideal(x^3+y^3+z^3); CF ideal(x*y*(x+y))
Fulton class : 3H
Fulton class : 3H
\end{verbatim}
while the Chern-Schwartz-MacPherson classes for the same ideals differ:
\begin{verbatim}
i15 : CSM ideal(x^3+y^3+z^3); CSM ideal(x*y*(x+y))
Chern-Schwartz-MacPherson class : 3H
                                    2
Chern-Schwartz-MacPherson class : 4H  + 3H
\end{verbatim}
The Euler characteristic in the second case is computed to be $4$, as it 
should. Taking the ideal $(xy(x+y)))$ in the {\em affine\/} plane 
gives a cone, so the Euler characteristic of the corresponding scheme in 
$\Abb^2$ must be $1$:
\begin{verbatim}
i17 : use QQ[x,y]; euleraffine ideal(x*y*(x+y))

o18 = 1
\end{verbatim}
while the Euler characteristic of the nonsingular affine cubic
$x^3+y^3=1$ is $-3$:
\begin{verbatim}
i19 : euleraffine ideal(x^3+y^3-1)

o19 = -3
\end{verbatim}
\end{example}

\begin{example}[A nonreduced example]
Here are $\cf$ and $\csm$ for a reduced pair of lines in $\Pbb^2$:
\begin{verbatim}
i20 : use ringP2; CF ideal(x*y); CSM ideal(x*y)
                 2
Fulton class : 2H  + 2H
                                    2
Chern-Schwartz-MacPherson class : 3H  + 2H
\end{verbatim}
In $\Pbb^3$, we can consider the ideal $(xy,xz,yz,z^2)=(x,z)(y,z)$:
this defines a scheme supported on two concurrent lines, but with a
nilpotent on the point of intersection. This can be checked with Macaulay2:
\begin{verbatim}
i23 : use ringP3; ass ideal(x*y,x*z,y*z,z^2)

o24 = {ideal (z, x), ideal (z, y), ideal (z, y, x)}

o24 : List
\end{verbatim}
And here are $\cf$ and $\csm$:
\begin{verbatim}
i25 : CF ideal(x*y,x*z,y*z,z^2); CSM ideal(x*y,x*z,y*z,z^2)
                 3     2
Fulton class : 4H  + 2H
                                    3     2
Chern-Schwartz-MacPherson class : 3H  + 2H
\end{verbatim}
As should be expected, $\cf$ detects the embedded component, while $\csm$ 
ignores it.
\end{example}

\begin{example}[Quintic threefold]
Fulton and Chern-Schwartz-MacPherson classes agree for nonsingular 
varieties $S$, as they both give the total (homology) Chern class of the 
tangent bundle of $S$. Here is the computation for the Fermat quintic in 
$\Pbb^4$:
\begin{verbatim}
i27 : use ringP4; quintic=ideal(x^5+y^5+z^5+w^5+t^5);

o28 : Ideal of ringP4

i29 : CF quintic; CSM quintic
                     4      3
Fulton class : - 200H  + 50H  + 5H
                                        4      3
Chern-Schwartz-MacPherson class : - 200H  + 50H  + 5H
\end{verbatim}
giving Euler characteristic$=-200$, as it should be. Computing the
Euler characteristic of singular quintic threefolds is equally
straightforward; here is a random example inspired by reading about
elliptic Calabi-Yau threefolds:
\begin{verbatim}
i31 : CSM ideal(x^3*t^2+x*z^4+w^5-y^2*t^3)
                                    4      3
Chern-Schwartz-MacPherson class : 4H  + 38H  + 5H
\end{verbatim}
that is, the hypersurface obtained by closing up $y^2=x^3+z^4x+w^5$ in
$\Pbb^4$ has Euler characteristic~4.
\end{example}

\begin{example}[Discriminants]
Identify $\Pbb^3$ with the space of triples of points in $\Pbb^1$. The
set of nonreduced triples forms a hypersurface of degree $4$. Here is its 
Chern-Schwartz-MacPherson class:
\begin{verbatim}
i32 : use ringP3; 

i33 : CSM ideal(-27*x^2*w^2+18*x*w*y*z+y^2*z^2-4*y^3*w-4*x*z^3)
                                    3     2
Chern-Schwartz-MacPherson class : 4H  + 6H  + 4H
\end{verbatim}
This agrees with the computation in \cite{MR99m:14103}. In general,
the Euler characteristic of the discriminant hypersurface for
$d$-tuples is $(d+1)$.

Identifying $\Pbb^5$ with the space of plane conics, we have similarly
a discriminant hypersurface parametrizing singular conics, that is,
pairs of lines; explicitly, this can be realized as the determinant of
a symmetric $3\times 3$ matrix. Its Chern-Schwartz-MacPherson class:
\begin{verbatim}
i34 : use QQ[x,y,z,w,t,u];

i35 : CSM ideal det matrix {{x,y,z},{y,w,t},{z,t,u}}
                                    5      4      3     2
Chern-Schwartz-MacPherson class : 6H  + 12H  + 14H  + 9H  + 3H
\end{verbatim}
The Chern-Schwartz-MacPherson class of the discriminant of plane {\em
  cubics\/} is computed in \cite{MR99m:14103}, Corollary~12; but
that computation seems to be computationally out of reach of {\tt
  CSM.m2} at present.
\end{example}

\begin{example}[$d$-tuples]
A good source of examples of applications of Segre classes is enumerative 
geometry. As the procedure described in \S\ref{segre} computes the
Segre class precisely by solving a number of enumerative problems, it
is hardly surprising that the enumerative answers can be decoded back
from the Segre class; the examples that follow illustrate this
procedure.

The degree of the PGL(2)-orbit closure of a configuration of $d$
points in $\Pbb^1$ (counting multiplicities) has been studied in
\cite{MR94j:14044}. For a configuration $C$, this degree computes the
number (with multiplicities) of translates of $C$ which contain three
given general points; the `predegree' of an orbit closure counts such
translates according to automorphisms of the $d$-tuple. In order to
use a Segre class to compute this predegree, one can parametrize
translates of a fixed $C$ by the $\Pbb^3$ of $2\times 2$ matrices
$$\pmatrix x & y \\ z & w \endpmatrix\quad;$$
the condition that the translate of $C$ contains a point determines a
surface in this $\Pbb^3$, and the predegree is given by the number of
points of intersection of three general such surfaces. The problem of
computing this number is not immediately reduced to B\'ezout's theorem
because these surfaces have an {\em excess intersection\/}. In
general, contributions of excess intersections can be evaluated in
terms of a Segre class by using Proposition~9.1.1 in
\cite{MR85k:14004}. With this in mind, the predegree of the orbit
closure is given by 
$$d^3-\int (1+dH)^3 i_* s(S,\Pbb^3)\quad,$$
where $d$ is the degree of $C$, and $S$ is the base scheme of the map 
$\Pbb^3 \dashrightarrow \Pbb^d$ mapping a matrix as above to the 
corresponding translate of $C$.

For a concrete example, consider the $5$-tuple with ideal generated by
$$s(s+3t)^2(s+5t)(s+16t)$$
in $\Pbb^1$. The reader should have no difficulties obtaining the
ideal of $S$. Using this, our routine computes $i_* s(S,\Pbb^3)$ as
$$\left\{\aligned
13H^2-70H^3 &\\
11H^2-58H^3 &\\
 9H^2-34H^3 &\\
 7H^2-22H^3
\endaligned\right.$$
in characteristic $2$, $3$, $5$, and $7$ respectively. For example
(the ideal is loaded from a separate file):
\begin{verbatim}
i36 : use ZZ/3[x,y,z,w]; load "dtupleideal.m2";
--loaded dtupleideal.m2

i38 : segre dtupleideal
                   3      2
Segre class : - 58H  + 11H
\end{verbatim}
Using the formula given above, these classes correspond to predegrees
$0$, $18$, $24$, $42$ respectively in char.~$2$, $3$, $5$, $7$. These
numbers are nicely explained by the result in \cite{MR94j:14044}: in
characteristic $2$ the tuple collapses to a pair of points, hence its
orbit closure has dimension~2; in characteristic $3$ it consists of
three points with multiplicities $3,1,1$; in characteristic $5$, three
points with multiplicities $2,2,1$; and in characteristic $7$ (and
most others, including~0) of four points, one of which double. These
multiplicities determine the predegree of the orbit closure, by
\cite{MR94j:14044}, Proposition~1.3; applying that result gives the
same predegrees as obtained here by brute force.
\end{example}

\begin{example}[Milnor classes]\label{milnor}
The function {\tt milnor} computes both Fulton and
Chern-Schwartz-MacPherson classes, giving $i_*$ of the difference
$$\csm(S)-\cf(S)\quad.$$
This class (up to a sign, cf.~the definition of $\mathcal M(Z)$ in
\cite{MR1795550}, p.~64) has been named the {\em Milnor class\/} of
$S$; to our knowledge, it has not been studied in any depth for
schemes other than reduced local complete intersection.

If $S$ is a hypersurface with isolated singularities, then $i_*$ of
the Milnor class of $S$ is simply (up to sign) $\mu H^n$, where $\mu$
is the sum of the Milnor numbers of the singularities; this is the
reason for the choice of terminology.
\begin{verbatim}
i39 : use ringP2; milnor ideal(y^6+z*x^3*y^2+z^2*x^4)
                    2
Fulton class : - 18H  + 6H
Chern-Schwartz-MacPherson class : 6H
                  2
Milnor class : 18H
\end{verbatim}
This says that the sum of the Milnor numbers of the curve
$y^6+x^3y^2z+x^4z^2=0$ in $\Pbb^2$ is~18. It may be checked that this
curve has singularities at $(x:y:z)=(1:0:0)$ and $(0:0:1)$, with
Milnor numbers respectively $3$ and $15$, consistently with this
information.

More generally, the coefficient of $H^n$ in the output of {\tt
  milnor\/} for an arbitrary hypersurface of $\Pbb^n$ computes Adam
  Parusi\'nski's generalization of the Milnor number,
  \cite{MR89k:32023}, whether the singularities are isolated or not. 

Our routine compute a notion of Milnor class for arbitrary projective
schemes. For example, the following would be the computation of the
Milnor class of the union of a line and a plane in $\Pbb^3$:
\begin{verbatim}
i41 : use ringP3; milnor ideal(x*y,x*z)
                 3     2
Fulton class : 2H  + 4H  + H
                                    3     2
Chern-Schwartz-MacPherson class : 4H  + 4H  + H
                 3
Milnor class : 2H
\end{verbatim}
In fact such computations may be performed in any characteristic; so
far as we know, no interpretation of the class is known in positive
characteristic.

Recent work of J\"org Sch\"urmann, \cite{math.AG/0202175}, relates
Milnor classes of complex local complete intersection with his
generalization of Deligne's functor of vanishing cycles.
\end{example}

\end{document}